\documentclass{svproc}

\usepackage{mathptmx}       % selects Times Roman as basic font
\usepackage{helvet}         % selects Helvetica as sans-serif font
\usepackage{courier}        % selects Courier as typewriter font
\usepackage{type1cm}        % activate if the above 3 fonts are
\usepackage{todonotes}
\usepackage{makeidx}         % allows index generation
\usepackage{graphicx}        % standard LaTeX graphics tool
\usepackage{multicol}        % used for the two-column index
\usepackage[bottom]{footmisc}% places footnotes at page bottom
\usepackage{amsmath, amsfonts, amssymb}
\usepackage{color}
\usepackage{pgf,tikz,pgfplots}
\usepackage[]{algorithm2e}
\usepackage{hyperref}

\usepackage{lineno}
\modulolinenumbers[5]
\linenumbers

\usepackage{subfig}

\makeindex             % used for the subject index
\newtheorem{Proposition}{Proposition}

\newtheorem{Remark}{Remark}

\newcommand{\dt}{\,\partial_t\, }
\newcommand{\dx}{\,\partial_x\, }
\newcommand{\dv}{\,\partial_v\, }

\newcommand{\dxx}{\,\partial_{xx}\, }

\newcommand{\dxxx}{\,\partial_{xxx}\, }

\newcommand{\dd}{\,\mathrm{d}}

\newcommand{\N}{\mathbb{N}}

\newcommand{\R}{\mathbb{R}}

\newcommand{\ipd}{{i+\frac{1}{2}}}
\newcommand{\imd}{{i-\frac{1}{2}}}

\newcommand{\M}{\mathcal{M}}
\newcommand{\T}{\mathcal{T}}

\newcommand{\ve}{\varepsilon}

\newcommand{\overbar}[1]{\mkern 1.5mu\overline{\mkern-1.5mu#1\mkern-1.5mu}\mkern 1.5mu}

\newcommand{\Thomas}[1]{\textcolor{black}{#1}}
\newcommand{\Oeps}[1]{\mathcal{O}(\varepsilon^{#1})}

\begin{document}
\mainmatter
\nolinenumbers
\title{Hybrid Kinetic/Fluid numerical method for the Vlasov-Poisson-BGK equation in the diffusive scaling}
\titlerunning{Hybrid Kinetic/Fluid numerical method}
\author{Tino Laidin \inst{1} \and Thomas Rey\inst{1}}

\institute{
Univ. Lille, CNRS, Inria, UMR 8524 - Laboratoire Paul Painlevé, F-59000 Lille, France
}

\maketitle

\abstract{
  This short note presents an extension of the hybrid, model-adaptation method introduced in [T.~Laidin, \textit{arXiv 2202.03696}, 2022] for linear collisional kinetic equations in a diffusive scaling to the nonlinear mean-field Vlasov-Poisson-BGK model. The aim of the approach is to reduce the computational cost by taking advantage of the lower dimensionality of the asymptotic model while reducing the overall error. It relies on two criteria motivated by a perturbative approach to obtain a dynamic domain adaptation. The performance of the method and the conservation of mass are illustrated through numerical examples.
\keywords{\Thomas{Vlasov-Poisson-BGK equations}; Diffusion scaling; Asymptotic preserving scheme; Micro-macro decomposition; Hybrid solver
%\\[5pt]
%{\bf MSC }(2010){Primary 65M08, 82M12; Secondary 35B40, 65M55} 
}
% msc 2010 classification of your contribution
% the list of the msc is available here: http://www.ams.org/mathscinet/msc/
}

\section{The model}
In this work, \Thomas{we are interested in the kinetic description of a system of particles interacting via both mean-field electromagnetic interaction and  collisions. Such system can be modeled using} the Vlasov-Poisson-BGK equation. The unknown is the probability distribution function $f=f(t,x,v)\in\R^+$ solution to:
\begin{equation}\label{VPBGK}
    \left\lbrace\begin{aligned}
        &\frac{\partial}{\partial t} f^\ve + \frac{v_x}{\varepsilon}\,\dx f^\ve + \frac{E^\ve}{\varepsilon} \,\partial_{v_x}f^\ve=\frac{1}{\varepsilon^2}\mathcal{Q}(f^\ve),\\
        &f(0,x,v)=f_0(x,v),
    \end{aligned}\right.
\end{equation}
where $(t,x,v)\in\R^+\times\Thomas{[0, x_*]}\times\R^{d_v}$ \Thomas{with periodic boundary conditions in the physical space} and $v = (v_x, v_y, v_z)$. The long range interactions are modeled via the self-consistent electrical field $\Thomas{E=E(t,x)}$ solution to the \Thomas{Poisson equation}
\begin{equation}\label{Poisson}
  \dx E^\ve=\rho^\ve-\overbar{\rho}\quad\text{with}\quad\overbar{\rho}=\iint f_0 \dd x\dd v.
\end{equation}
The short-range interactions between particles are taken into account through the linear BGK\Thomas{-like} relaxation collision operator
\begin{equation*}
    \mathcal{Q}(f^\varepsilon)(t,x,v)=\rho^\ve(t,x)\M(v)-f^\varepsilon(t,x,v)\Thomas{, \quad \forall (t,x,v)\in\R^+\times\Thomas{[0, x_*]}\times\R^{d_v},}
\end{equation*}
where \Thomas{the global Maxwellian and local density  are respectively defined as}
\begin{equation*}
    \M(v)=\frac{e^{-|v|^2/2}}{(2\pi)^{d_v/2}}, \quad \rho^\varepsilon(t,x)=\int f^\varepsilon(t,x,v)\dd v =: \langle f^\ve(t,x,\cdot)\rangle.
\end{equation*}
In \eqref{VPBGK}, the scaling parameter $\varepsilon$ is related to the Knudsen number: the ratio between the mean free path of the particles and the length scale of observation. It is now well known \cite{GoudonPoupaud2001} that the limit case $\varepsilon=0$ is described by a drift-diffusion equation on the density $\rho$: \Thomas{when $\ve \to 0$, the distribution function $f^\ve$ converges towards a Maxwellian distribution $\rho \M$ whose density $\rho$ is solution to} %is: $f^\varepsilon\underset{\varepsilon\rightarrow0}{\longrightarrow}\rho\M$ and
\begin{equation*}\label{DD}\tag{$P$}
    \left\lbrace\begin{aligned}
        &\dt \rho -\dx J = 0,\quad J=\dx\rho - E\rho,\\
        &\rho(0,x)=\rho_0(x)\Thomas{, \quad \forall x \in [0, x_*]}.
    \end{aligned}\right.
\end{equation*}
The approximation of solutions to \eqref{VPBGK}-\eqref{Poisson} can be \Thomas{computationally expensive because of the high dimensionality of the kinetic phase space. Nevertheless, using the more accurate kinetic description of the system may not be necessary in the whole computing domain because of the local validity of the fluid description (the system being in a state of thermodynamic equilibrium). This latter is less accurate in describing the kinetic zones, but computationally  less costly.} The aim of this work is therefore to design a hybrid kinetic/fluid scheme with an automatic domain adaptation method. It relies on a robust numerical scheme for the kinetic equation, on relevant criteria to carefully determine fluid and kinetic regions and on a smart implementation.

\section{Macroscopic models}\label{sec:MacroModel}
The aim of this section is to derive a higher-order macroscopic model from which we deduce a macroscopic coupling criterion. \Thomas{It generalizes the approach presented in \cite{Laidin2022}.}
Let us now introduce the truncated Chapman-Enskog expansion of the distribution function $f^\varepsilon$ at order $K\in\N^*$: 
\begin{equation}\label{eq:CE}
    f^\varepsilon(t,x,v)=\rho^\varepsilon(t,x)\M(v) + \sum^K_{k=1}\varepsilon^k h^{(k)}(t,x,v).
\end{equation}
By inserting \eqref{eq:CE} into the original equation \eqref{VPBGK}, one can identify powers of epsilon to obtain 
\begin{subequations}\label{Identification}
\begin{alignat}{2}
k=0:&\quad &h^{(1)} = &-\T(\rho^\varepsilon\M),\\
k=1:&\quad &h^{(2)} = &-\dt(\rho^\varepsilon\M) -\T(h^{(1)}),\\
2\leq k\leq K-1:&\quad &h^{(k+1)} = &-\dt h^{(k-1)} - \T(h^{(k)}),
\end{alignat}
\end{subequations}
where $\mathcal{T}f=v_x \, \partial_xf+E \,\partial_{v_x} f$ is the transport operator. 
To obtain a hierarchy of macroscopic models, one considers different truncation orders $K$, then plugs \eqref{eq:CE} into \eqref{VPBGK} and integrates in velocity. The order $K=1$ allows us to (formally) recover the asymptotic model \eqref{DD}. 

Let us now present the idea behind the computations in a $1D/1D$ setting for the case $K=3$. Note that the same method can be applied up to the full $3D/3D$ setting. \Thomas{Integrating in velocity then yields}
\begin{equation}\label{VPBGKintegrated}
\dt\rho\Thomas{^\ve}+\dx\langle vh^{(1)}\rangle+\varepsilon^2\dx\langle vh^{(3)}\rangle=\Oeps{4}.
\end{equation}
 The set of equations \eqref{Identification} allows us to compute the functions $h^{(k)}$, $k=1,2,3$. Using the identity $\dv\M=-v\M$, \Thomas{one has}
\begin{equation*}
    h^{(1)}=-v\M J^\varepsilon,\quad J^\ve=\dx\,\rho^\ve-E^\ve\rho^\ve.
\end{equation*}
\begin{equation*}
    h^{(2)}=-\M\dt\rho^\varepsilon+v^2\M\dx J^\varepsilon+(1-v^2)\M E^\varepsilon J^\varepsilon,
\end{equation*}
and
\begin{equation*}
\begin{aligned}
    h^{(3)}=&v\M J^\varepsilon+v\M\dx(\dt\rho^\varepsilon)-v^3\M\dxx J^\varepsilon-(v-v^3)\M\dx(E^\varepsilon J)\\
    &-v\M E^\varepsilon\dt\rho^\varepsilon -(2v-v^3)\M E^\varepsilon\dx J^\varepsilon -(v^3-3v)(E^\varepsilon)^2J^\varepsilon.
\end{aligned}
\end{equation*}
It remains to compute the quantities $\dx\langle vh^{(1)}\rangle$ and $\dx\langle vh^{(3)}\rangle$. By the definition of $\M$, one can explicitly compute its second and fourth moments: $m_2=1$ and $m_4=3$. It yields
\begin{equation*}
    \dx\langle vh^{(1)}\rangle=\dx J^\varepsilon.
\end{equation*}
Moreover, to avoid some approximation of mixed derivatives, we observe from \eqref{VPBGKintegrated} that 
\begin{equation}\label{dtRho}
    \dt\rho^\ve=\dx J^\ve + \Oeps{2},
\end{equation}
and we deduce that
\begin{equation}\label{dtJ}
    \dt J^\ve=\dxx J^\ve-E^\ve\dx J^\ve-\rho^\ve\dt E^\ve+\Oeps{2}.
\end{equation}
Replacing the time derivatives by their approximations \eqref{dtRho} and \eqref{dtJ} yields
\begin{equation*}
\begin{aligned}
    \dx\langle vh^{(3)}\rangle=&-\dxxx J^\varepsilon+\dx\rho^\varepsilon\dt E^\varepsilon-\rho^\varepsilon\dt(\dx E^\varepsilon)+3\dx E^\varepsilon\dx J^\varepsilon \\
    &+E^\varepsilon\dxx J^\varepsilon+2J^\varepsilon\dx(\dx E)+\Oeps{2}.
\end{aligned}
\end{equation*}
Finally, using the Poisson equation \eqref{Poisson} and rearranging the terms, we obtain a higher order macroscopic model:
\begin{equation}\label{HigherOrderModel}
    \dt\rho^\varepsilon-\dx J^\varepsilon=-\varepsilon^2\mathcal{R}+\Oeps{4},
\end{equation}
where the remainder $\mathcal{R^\varepsilon}$ is given by
\begin{equation}\label{Remainder}
    \mathcal{R}^\varepsilon=-\dxxx J^\varepsilon + E\dxx J^\varepsilon + (2\rho^\varepsilon-3\overbar{\rho})\dx J^\varepsilon +2J^\varepsilon\dx\rho^\varepsilon -\dx\rho^\varepsilon\dt E^\varepsilon.
\end{equation}
Furthermore, we emphasize the fact that this remainder term does not depend on the velocity variable but on the macroscopic quantities $\rho$ and $E$. It quantifies very accurately the deviation from the thermodynamical equilibrium.

\section{Numerical scheme}\label{sec:Scheme}
In this part, we briefly recall the derivation of the micro-macro model for \eqref{VPBGK} introduced in \cite{CrouseillesLemou2011}. Then, we present a micro-macro finite volume scheme that enjoys the property of being Asymptotic Preserving \Thomas{(AP)}, which is a crucial point of the hybrid method we present. 

Let us decompose the distribution $f^\varepsilon$ as follows:
\begin{equation}\label{rhoM+g}
    f^\varepsilon = \rho^\varepsilon\M+g^\varepsilon.
\end{equation}
We introduce the orthogonal projector $\Pi$ in $L^2(\dd x\dd v\M^{-1})$ on $Ker(\mathcal{Q})$ defined for all $f\in L^2(\dd x\dd v\M^{-1})$ by: $$\Pi f=\langle f\rangle\M.$$
\Thomas{The so-called} micro equation is obtained by plugging \eqref{rhoM+g} into \eqref{VPBGK} and applying $(I-\Pi)$. \Thomas{Moreover}, plugging \eqref{rhoM+g} into \eqref{VPBGK} and applying the projection $\Pi$ yields the macro equation. The micro-macro model is then given by:
\begin{gather}\label{Micro}\tag{$Micro$}
    \dt g^\varepsilon+\frac{1}{\varepsilon}(\T g^\varepsilon - \dx\langle v_xg^\varepsilon\rangle\M + v_x\M J^\varepsilon) =\frac{-1}{\varepsilon^2}g^\varepsilon, \\
\label{Macro}\tag{$Macro$}
    \dt\rho^\varepsilon+\frac{1}{\varepsilon}\dx\langle v_xg^\varepsilon\rangle=0.
\end{gather}

Let us now present the discretization of the \eqref{Micro}-\eqref{Macro} model. We shall adopt a finite volume approach to discretize the phase space and \Thomas{present only the $1D/3D$ case that we shall consider in our numerical simulations}.% problem, namely one periodic dimension in position and three in velocity. In particular, let $\Omega_x=[0,x_\star]$. 

% In this setting, the \eqref{Micro}-\eqref{Macro} equations reduces to 
% \begin{equation*}
% \left\lbrace\begin{aligned}
%     & \dt g^\varepsilon + \frac{1}{\varepsilon^2}g^\varepsilon + \frac{1}{\varepsilon}\big(\T g^\varepsilon - \dx\langle v_xg^\varepsilon\rangle\M + v_x\M J^\varepsilon_x\big)=0\\
%     &\dt\rho^\varepsilon + \frac{1}{\varepsilon}\dx\langle v_xg^\varepsilon\rangle=0,
% \end{aligned}\right.
% \end{equation*}
\paragraph{The mesh.} The velocity domain is restricted to the bounded symmetric cube $[-v_\star,v_\star]^3$. We consider a Cartesian mesh of the phase space composed of $N_v=2L$ velocity cells in each direction arranged symmetrically around $v=0$. Let $\mathcal{J}=\lbrace -L+1,\dots,L \rbrace$ and  $j=(j_x,j_y,j_z)\in\mathcal{J}^3$ be a multi-index. The cells of the velocity mesh are denoted by $\mathcal{V}_j$ for $j\in \mathcal{J}^3$. Each cell $\mathcal{V}_j$ has a constant volume $\Delta v^3$ and midpoint $v_j$.

\Thomas{The physical domain, a torus $\mathbb{T}$ of length $x_*$, is discretized into} $N_x$ primal cells $$\mathcal{X}_i=(x_\imd,x_\ipd), \quad i\in\mathcal{I}=\mathbb{Z}/N_x\mathbb{Z},$$ of constant length $\Delta x$ and centers $x_i$. We also define dual cells $\mathcal{X}_\ipd=(x_i,x_{i+1})$ for $i\in\mathcal{I}$, of constant length $\Delta x$ and centers $x_\ipd$. \Thomas{This defines respectively the so-called primal and staggered meshes as}
$$K_{ij}= \mathcal{X}_i\times\mathcal{V}_j\text{ and }K_{\ipd,j}= \mathcal{X}_\ipd\times\mathcal{V}_j,\quad\forall(i,j)\in\mathcal{I}\times\mathcal{J}^3.$$
Finally, we set a time step $\Delta t>0$ and define $t^n=n\Delta t$ for $n\in\mathbb{N}$.
The Maxwellian is discretized as a product of 1-dimensional Gaussians in such a way that it satisfies \Thomas{discrete counterparts of its continuous properties}, namely parity, positivity, and unit mass.

The transport terms in \eqref{Micro} are approximated using a first order upwind scheme and the time derivative is dealt with using a first order exponential time integrator. Regarding the macro equation, we take advantage of approximating the perturbation on the dual cells in position for the space derivative. The usual choice \cite{BennouneLemouMieussens2008} to obtain an AP scheme is to implicit the stiff term of \eqref{Macro}. The equation on the electrical field is solved using a centered finite difference scheme. The numerical scheme reads as follows:
\begin{Proposition}
    Let $n\in \mathbb{N}$. Let $\left(g_{\ipd,j}^{\varepsilon,n}\right)_{ij}$ and $\left(\rho^{\varepsilon,n}_i\right)_{i}$ be given by the following micro-macro finite volume scheme:
    \begin{gather*}
        g_{\ipd,j}^{\varepsilon,n+1} = g_{\ipd,j}^{\varepsilon,n}e^{-\Delta t/ \varepsilon^2} - \varepsilon(1-e^{-\Delta t/ \varepsilon^2})\left(\frac{T^{\varepsilon,n}_{\ipd,j}}{\Delta x\Delta v^3}+\xi_j\M_jJ^{\varepsilon,n}_\ipd\right), \\
        \rho_i^{\varepsilon,n+1} = \rho_i^{\varepsilon,n} - \frac{\Delta t}{\varepsilon\Delta x}\left(\langle \xi g^{\varepsilon,n+1}_\ipd \rangle_\Delta-\langle \xi g^{\varepsilon,n+1}_\imd \rangle_\Delta \right),\\
        E^{\ve,n}_\ipd-E^{\ve,n}_\imd = (\rho_i^{\ve,n}-\overbar{\rho})\Delta x,
    \end{gather*}
    where $\xi_j=\xi_{(j_x,j_y,j_z)}=v_{j_x}\,\forall j\in\mathcal{J}^3$ and $T^{\varepsilon,n}_{\ipd,j}$ is the discretization of the transport terms.
    Assuming some uniform bounds in $\varepsilon$ on $\rho^\varepsilon$ and for a fixed mesh size $\Delta x$, $\Delta v>0$, the scheme enjoys the AP property in the diffusion limit. This property does not depend on the initial data, and the associated limit scheme reads
    \begin{gather*}
        \rho_i^{n+1} = \rho_i^{n} +  m_2^{\Delta v}\frac{\Delta t}{\Delta x}\left(J_\ipd^{n}-J_\imd^{n}\right),\\
        E^{n}_\ipd-E^{n}_\imd = (\rho_i^{n}-\overbar{\rho})\Delta x,
    \end{gather*}
    with the limit flux
    \begin{gather*}
        J_\ipd^{n} = \frac{\rho^n_{i+1}-\rho^n_{i}}{\Delta x} - E_\ipd\rho^n_\ipd,
    \end{gather*}
    where $m_2^{\Delta v}$ is given by $\sum_{l\in\mathcal{J}}v_l^2M_l\Delta v_l$.
\end{Proposition}

\section{Hybrid method}
Let us now introduce \Thomas{the main contribution of this work, namely} the hybrid method between kinetic and fluid schemes. \Thomas{It consists in} a coupled solver that is faster than a full kinetic one to solve \eqref{VPBGK} while still being accurate. %These methods come naturally when designing efficient numerical codes while ensuring reasonable computation times.
Following \cite{FilbetRey2015,Laidin2022} we first construct a hybrid kinetic/fluid solver with a dynamic domain adaptation method and present its implementation.

\subsection{Coupling criteria}
The idea of the dynamic domain adaptation method is twofold. On the one hand, the subdomains must accurately describe the state of the solution. In particular, the fluid model has to be used only where the solution is near the local velocity equilibrium. On the other hand, the method has to be dynamic, in the sense that the subdomains are adapted at each time step. To determine in which domain each cell lies, we use two criteria based on the higher order macroscopic model \eqref{HigherOrderModel} introduced in Section \ref{sec:MacroModel} and the norm of the perturbation $g^\varepsilon=f^\varepsilon-\rho^\varepsilon\M$. Let $\mathcal{R}_i^{\varepsilon,n}$ be a discretization of the remainder \eqref{Remainder} at time $t^n$ in cell $\mathcal{X}_i$. 

The coupling procedure will unfold as follows: $\mathcal{R}_i^{\varepsilon,n}$ is computed using both the kinetic density $\rho^\varepsilon$ in kinetic cells and the fluid density $\rho$ in fluid cells. Cell changes happen when one of the following situations occurs, depending on a coupling threshold investigated in the previous work \cite{Laidin2022}:
\begin{itemize}
    \item If $\mathcal{R}_i^{\varepsilon,n}$ is small (w.r.t. a fixed threshold), the solution \eqref{HigherOrderModel} is close to the limit model \eqref{DD} and one must use it;
    \item If $g^\varepsilon$ is small (w.r.t. another fixed threshold), the solution is close to the local equilibrium in velocity, and one must consider the fluid description.
\end{itemize}

%Let $\eta_0,\,\delta_0>0$ be coupling thresholds. For a cell $\mathcal{X}_i$ at time $t^n$, we consider that whenever both $\left|\mathcal{R}^{\varepsilon,n}_i\right|$ and $ \|g^{\varepsilon,n}_{\ipd}\|_{2}$ are small enough with respect to the coupling thresholds, it can remain a fluid cell or become one. 
If any of those two conditions are not met then, depending on its current state, the cell either stays or becomes kinetic.

\begin{Remark}
Note that in a kinetic cell, the criterion on the norm of $g^\varepsilon$ is mandatory. Indeed, the remainder $\mathcal{R}^{\varepsilon,n}_i$ could be small because of small gradients, but the perturbation large. In this situation, one does not want to change from kinetic to fluid. As an example, one could take a distribution function at equilibrium in position and far from the Maxwellian in velocity (see \cite{FilbetRey2015} for details).
\end{Remark}

\subsection{Implementation}
The way the method is implemented is crucial. Indeed, the main goal of the method is to avoid the update of the perturbation which is the most computationally expensive part of the code. Therefore, from an implementation point of view, $g^\varepsilon$ is not updated in fluid regions and is set to $0$ only when needed. In particular, it occurs only when a fluid cell becomes kinetic and when saving data. The interface conditions between kinetic and fluid cells are treated in the same way as in \cite{Laidin2022}.

\begin{Remark}
    Note that we want to start the resolution with the approach containing the full information on the system. Hence, the domain is initialized as fully kinetic. 
    Moreover, let us emphasize that in practice the kinetic fluxes are explicitly computed. Therefore, the hybrid setting is an explicit method.
\end{Remark}

\section{Numerical results}

Let us now present some numerical simulations with our approach. 

\begin{figure}
    \centering
    \includegraphics[width=0.65\textwidth]{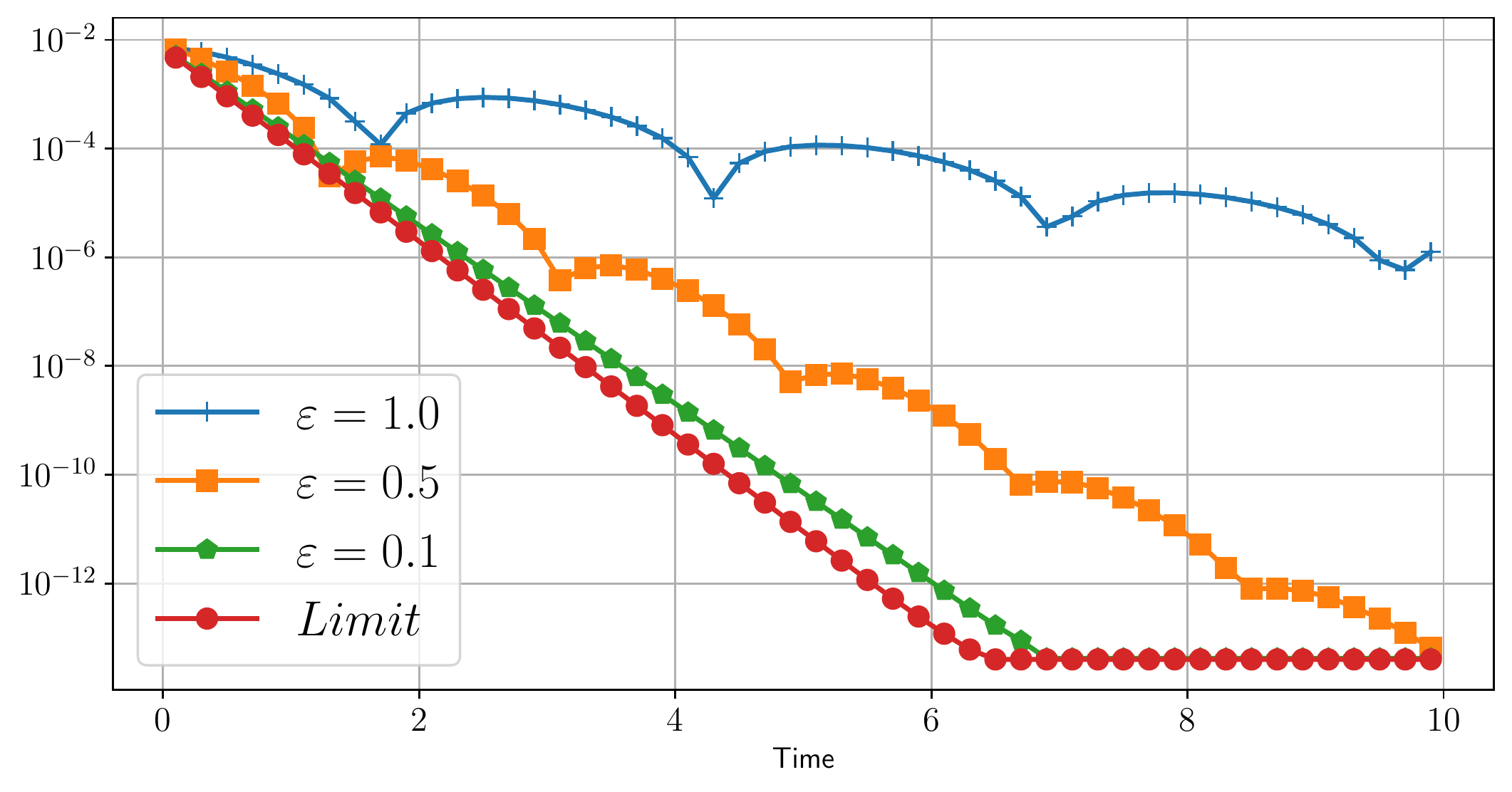}
    \caption{\textbf{Fully kinetic scheme: Landau damping.} Time evolution of $||E(t)||_2$ for different values of $\varepsilon$.}
    \label{fig:LandauEpsilon}
\end{figure}

We shall start with the validation of the AP property of the micro-macro Vlasov-Poisson-BGK solver with exponential integrator presented in Section \ref{sec:Scheme}. It is a combination of the methods presented in \cite{Lemou2010,CrouseillesLemou2011} that have never been implemented in any work, to the best of our knowledge.
Figure \ref{fig:LandauEpsilon} presents the time evolution of the $L^2$ norm of the energy of solutions to equations (\ref{VPBGK}--\ref{Poisson}) with different values of the relaxation parameter $\ve$, for the seminal weak Landau damping initial data from \cite{CrouseillesLemou2011}.
We observe the convergence with respect to $\ve$ of the energy. The oscillations due to the Vlasov-Poisson transport term $\T$ occur only in the kinetic regime, when $\ve$ is large. They are then damped by the linear BGK term for smaller values of $\ve$, where exponential decay of the electric field occurs.

\begin{figure}
    \centering
    \includegraphics[width=0.95\textwidth]{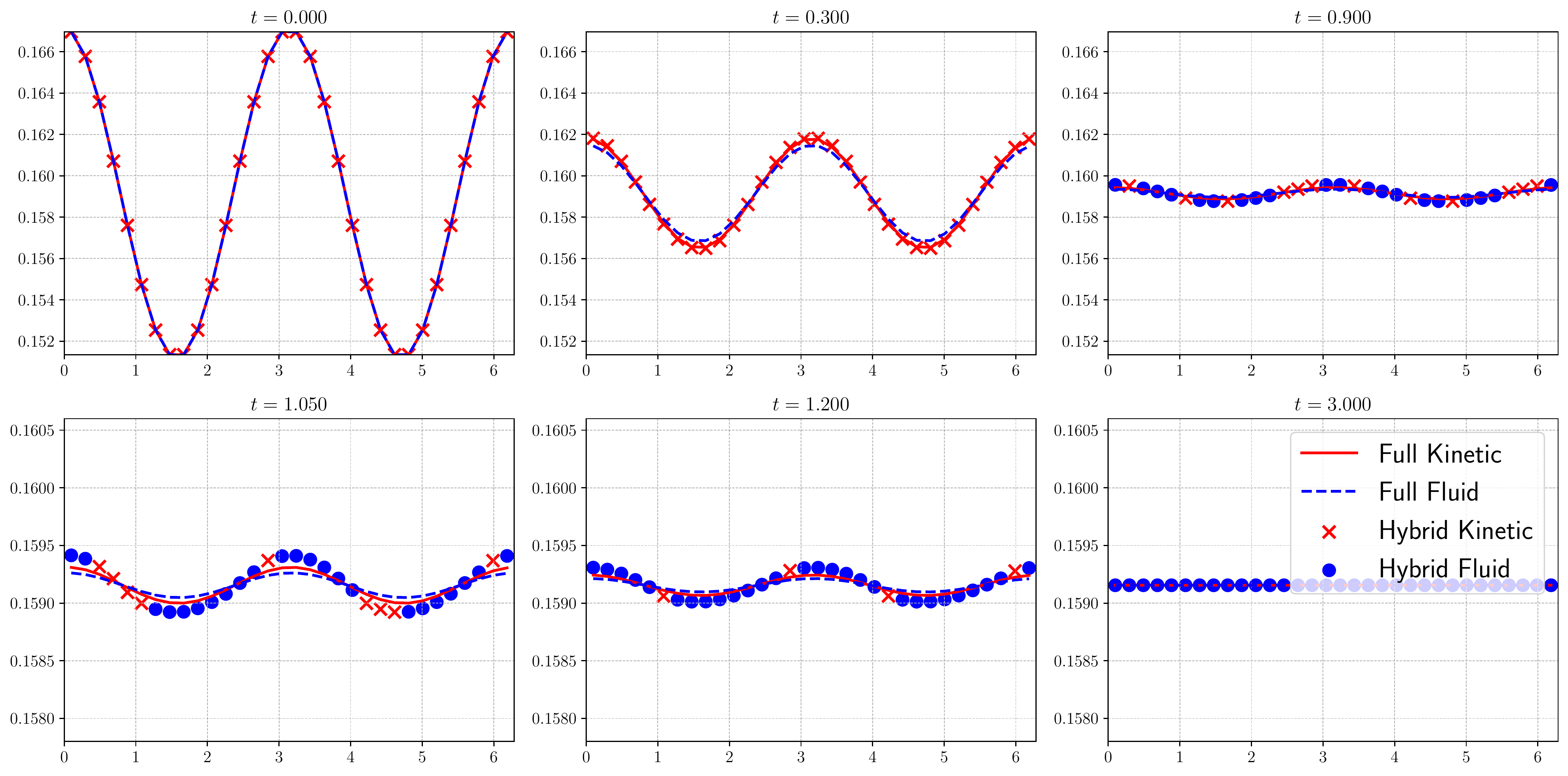}
    \caption{\textbf{Comparison of the solvers.} Time evolution of the space densities for $\varepsilon=0.1$ with a smooth initial data.}
    \label{fig:HybridPoissonComp01}
\end{figure}

We now turn our attention to the hybrid method. Figure \ref{fig:HybridPoissonComp01} assesses the validity of this new method by computing the time evolution of the density of a smooth solution. The initial condition is a Maxwellian distribution $\rho \M$, whose initial space dependent density is given by $\rho(x) = 1 + 0.05\cos(2x)$ for $x_* = 2\pi$ with a fixed $\ve = 0.1$. 
Although the solution is far from the fluid description because of this mild value of $\ve$, we observe an almost perfect agreement between the fully kinetic and the hybrid kinetic-fluid solvers.
We also observe the back-and-forth phenomenon between kinetic and fluid cells, resulting in large time in a full fluid (albeit accurate with the kinetic equation) solver for the global equilibrium. 
The speed-up tables are similar to those presented in \cite{Laidin2022}, where factors up to $~400$ have been observed between the hybrid and the fully kinetic solver.

%The speed up of this domain adaptation is presented in Table  \ref{tab:Speedup}. We observe the huge gain of the hybrid method compared to the fully kinetic one, with factors of the order 400. 
%\begin{table}
%    \centering
%    \begin{tabular}{|c|c|c|c|c|c|c|c|c|}
%        \hline
%        $\varepsilon$ & $10^{-4}$ & $10^{-3}$ & $5\times10^{-3}$ & $10^{-2}$ & $5\times10^{-2}$ & $10^{-1}$ & $5\times10^{-1}$ & $1.0$ \\
%        \hline
%        \,Speedup\, & $479.64$ & $464.12$ & $488.18$ & $13.91$ & $1.46$ & $1.02$ & $1.00$ & $0.97$\\
%        \hline
%    \end{tabular}
%    \caption{Speedup of the Hybrid method, initial data close to equilibrium.}
%    \label{tab:Speedup}
%\end{table}

\begin{figure}[ht!]
    \centering
    \includegraphics[width=0.95\textwidth]{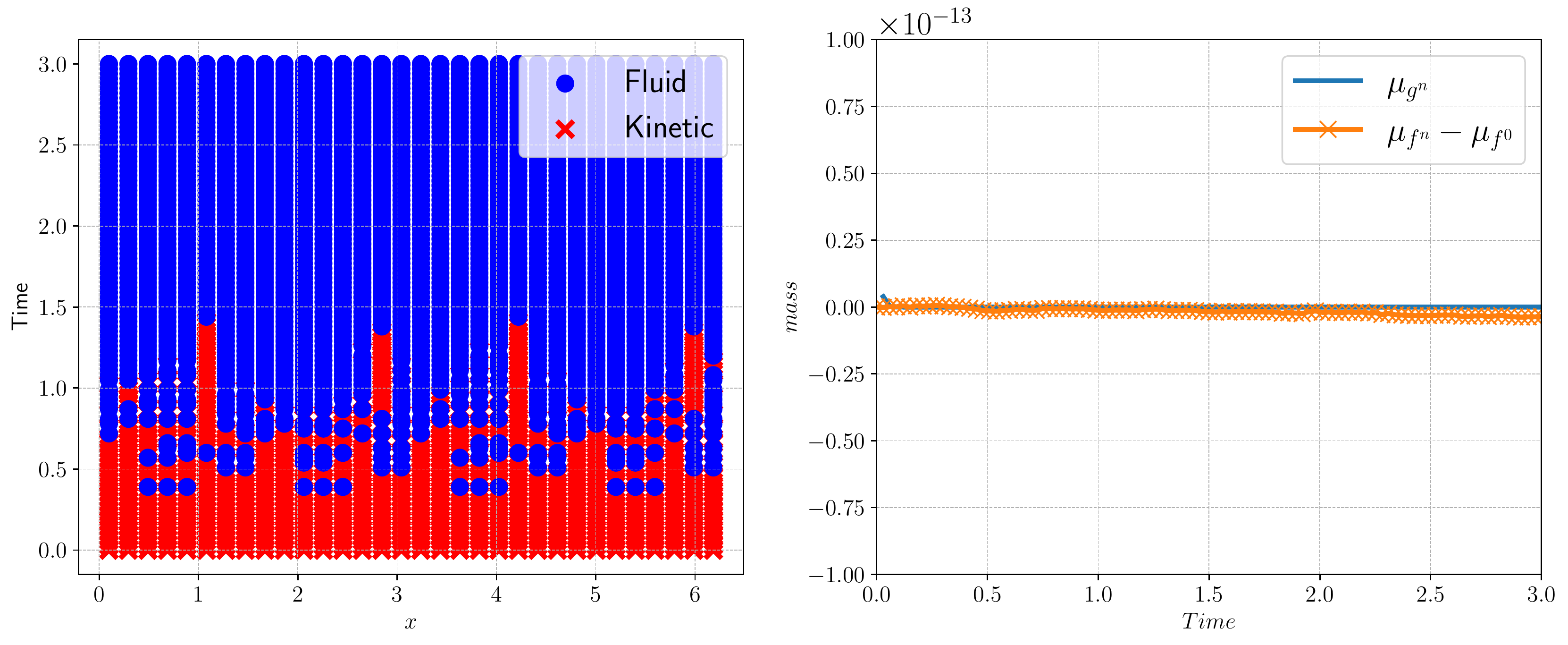}
    \caption{\textbf{Hybrid scheme.} \textit{Left.} Time evolution of the state of the cells. \textit{Right.} Evolution of the  mass variation (orange crosses) and mass of $g^\ve$ (solid blue line) for $\ve=0.1$.}
    \label{fig:State_MassVar_Poisson01}
\end{figure}

Then this domain adaptation is investigated in Figure \ref{fig:State_MassVar_Poisson01}. One can observe the precise domain adaptation during time. We notice in particular the quick vanishing of the kinetic cells in favor of the fluid ones (and hence a computational speedup). This adaptation phenomenon can bring mass variation as noticed in \cite{Laidin2022}, but we observe that it remains very close to the machine precision.

\begin{figure}
    \centering
    \includegraphics[width=0.65\textwidth]{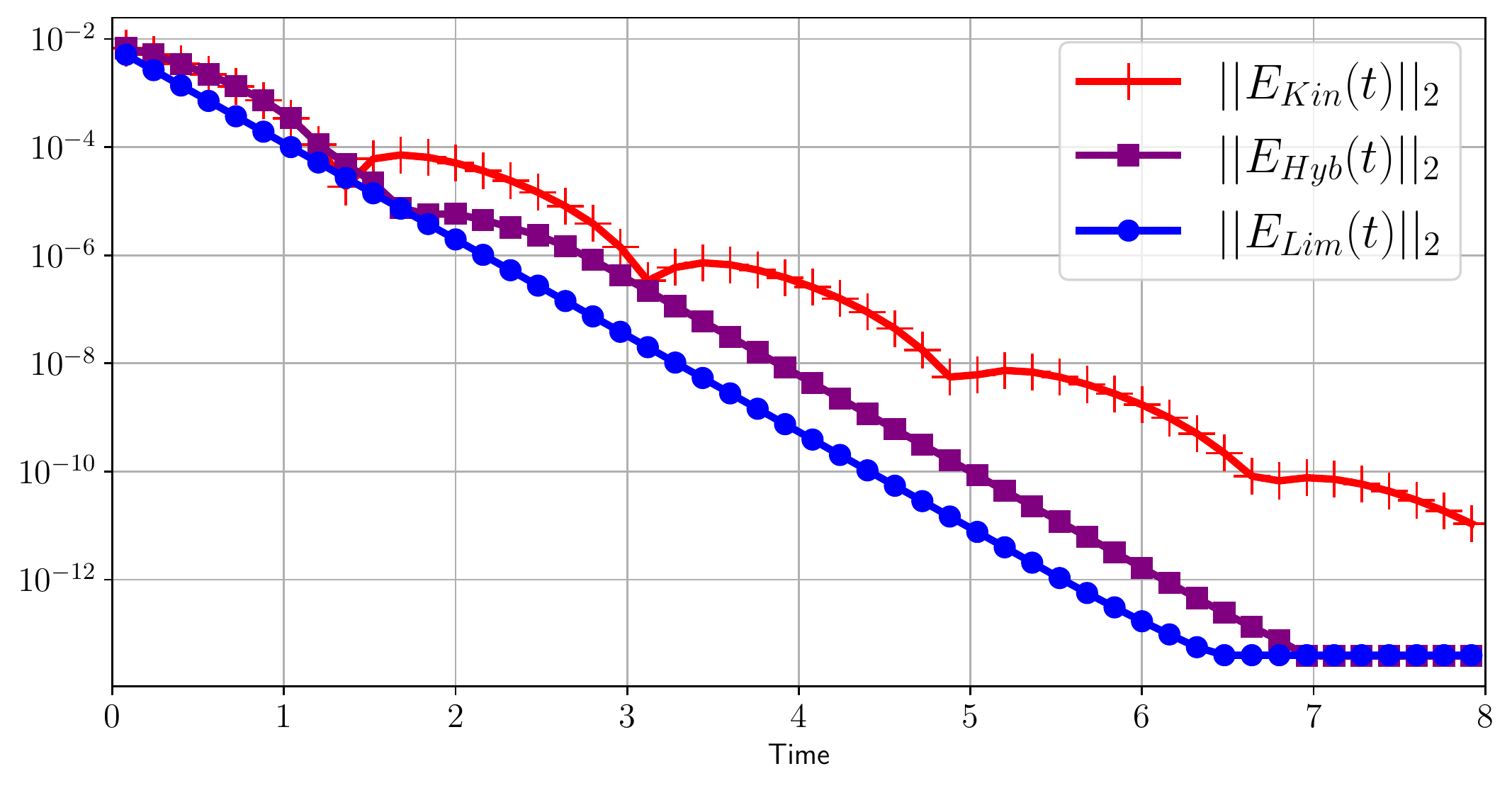}
    \caption{\textbf{Landau damping revisited.} Time Evolution of $||E(t)||_2$ for $\ve=0.5$.}
    \label{fig:LandauHybrid05}
\end{figure}

Finally, we investigate again the case of the weak Landau damping using the new hybrid solver. We observe in Figure \ref{fig:LandauHybrid05} that this new method is able to accurately capture the oscillations induced by the transport operator $\T$ in short time. Nevertheless, these oscillations are destroyed by the switch to a full fluid solver, which relax exponentially because of its relaxation structure.

\bigskip

\noindent
{\bf Acknowledgements.} 
The authors were partially funded by Labex CEMPI (ANR-11-LABX-0007-01) and the MSCA DN-2022 program DATAHYKING. 

% ===============================
% 	References
% ===============================

% the bibtex style to be used is spmpsci.bst
\bibliography{HybridMicroMacro2021.bib}
\bibliographystyle{spmpsci}

\end{document}